\documentclass[12pt]{article}
\usepackage[cp1251]{inputenc}
\usepackage[T2A]{fontenc}

\evensidemargin\oddsidemargin \topmargin2cm

\usepackage{amsfonts,amscd,amsmath}

\usepackage{
amssymb,amsmath,
latexsym,euscript
}

\usepackage{amsfonts,amscd,amsmath}

\def\d {\partial }

\begin{document}

\title{\textbf{A universal boundary value problem for\\
partial
differential equations}}

\author{V.Zh. Sakbaev $^1$\footnote{fumi2003@mail.ru} and I.V. Volovich
$^2$\footnote{volovich@mi.ras.ru},
 \\[2.7mm]
 ${}^1$\small{Moscow Institute of Physics and Technology
}\\
\small{Institutskiy per. 9, 141700, Dolgoprudny, Moscow reg., Russia}\\
${}^2$\small{Steklov Mathematical Institute}
\\
\small{Gubkin St. 8, 119991, Moscow, Russia}\\
}
\date{ }
\maketitle

\begin{abstract} A new boundary value problem for partial differential equations is discussed.  We consider an arbitrary solution of an elliptic or parabolic  equation in a given domain and no boundary conditions are assumed. We study which restrictions  the boundary values  of the solution
and its normal derivatives must satisfy.  Linear integral equations for the boundary values of the solution 
and its normal derivatives are obtained,  which we  call the universal boundary value equations. A universal boundary value problem is  defined as a partial differential equation together with  the  boundary data which specify the values of the solution on the boundary and its normal derivatives and satisfy  to the universal boundary value equations.

For the equations of mathematical physics such as  Laplace's and the heat equation the solution of the universal  boundary value problem is presented. Applications to cosmology and quantum mechanics are mentioned.

\end{abstract}

\newpage

\section {Introduction}

An ordinary boundary value problem is  defined as a partial differential equation in a domain together with  the boundary conditions.
A solution to the  boundary value problem is a solution to the differential equation which also satisfies the boundary conditions.

The form of the boundary conditions depends on the type of the partial differential equation (elliptic, parabolic or hyperbolic), see for example \cite{Vladimirov}.  Examples of boundary value
problems for  Laplace's equation are the Dirichlet problem, the Neumann problem and the third  boundary
value problem.

In this note a new  boundary value problem is discussed.
We consider an arbitrary   solution of an elliptic or parabolic equation in a given domain and no boundary conditions are assumed. We study which restrictions  the boundary values  of the solution
and its normal derivatives must satisfy.  Linear integral equations for the boundary values of the solution 
and its normal derivatives are obtained   which we  shall call the {\it universal boundary value equations}.

Our aim is to obtain a relation for the boundary
values of a function and its derivatives if the function is
a solution of the differential equation in the given domain.

If $u$ is a solution of  Laplace's equation $\Delta u=0$  in a domain $G$ with the boundary $S$ let $\gamma_0 u=u|_S$ be its trace on the boundary and $\gamma_1 u= \partial u/\partial \nu |_S$ the trace of the normal derivative. Then the universal boundary value equations have the form
$$
{\bf A}\gamma_0 u+{\bf B}\gamma_1 u=0,\,\,\,\, \int\limits_S\gamma_1 uds=0,
$$
where ${\bf A}$ and ${\bf B}$ are integral operators on $S$ (see Eqs (\ref{(2)}) and (\ref{(3)}) below). Conversely,
if $u_0$ and $u_1$ are two functions on $S$ which satisfy the universal boundary value equations
${\bf A} u_0+{\bf B} u_1=0, \,\int\limits_Su_1ds=0$ then there exists a solution $u$ of Laplace's equation $\Delta u=0$ such that $u|_S=u_0, \, \partial u/\partial \nu |_S=u_1$.

A very simple example  is given by the 1-dimensional "Laplace's equation"  \, \,\,\, $u^{\prime\prime}(x)=0 $, $~~~$ $a<x<b$. In this case the universal boundary value equations are the following two relations: $ u(a)-u(b)=u^{\prime}(a)(a-b) , \, u^{\prime}(a)=u^{\prime}(b) $ among four numbers (boundary values)
$u(a),u(b),u^{\prime}(a),u^{\prime}(b)$ (compare Eqs (\ref{(2)}) and (\ref{(3)}) below) .

For the equations of mathematical physics such as  Laplace's and the heat equation the solution of the universal  boundary value problem is presented. Applications to cosmology and quantum mechanics are mentioned.

\section{Universal boundary value problem for Laplace's equation}

Let us consider  Laplace's equation
\begin{equation}\label{(1)}
\Delta u=0,\ \, x\in G\subset R^3, \end{equation}
in the bounded domain $G$ of the Euclidean space $R^3$ with the boundary $S=\d G \in C^3$ which satisfies the additional condition of strong convexity type (an exact formulation of this condition is given below).

To study the boundary value problems for  Laplace's equation (1) we introduce the functional space
$C^{1,norm}(\bar G)$ which is the space of continuously differential functions on domain $G$ with the "proper normal derivative" ${{\d u}\over {\d \nu }}$ on the boundary $\d G$ (see \cite{Vladimirov}). Here $\nu $ is the continuous vector field of the unit exterior normal vectors  on the surface $S=\d G$.

{\bf Theorem 1}. {\it Let the function $u\in C^2(G)\bigcap
C^{1,norm}(\bar G) $ satisfies the equation (1) and the  surface $S\in C^3$ satisfies the strong convexity condition (see below). Then the boundary values $u_0=u|_S$ and
$u_1={{\d u}\over {\d \nu }}|_S$ satisfy the equations (universal boundary value equations)}
\begin{equation}\label{(2)}
u_0(x)+{1\over {2\pi }}\int\limits_S{{cos\phi _{xy}}\over {|x-y|^2}}u_0(y)ds_y-{1\over {2\pi }}\int\limits_S{1
\over {|x-y|}}u_1(y)ds_y=0,\  x\in S;
\end{equation}
\begin{equation}\label{(3)}
\int\limits_Su_1ds=0.
\end{equation}
{\it Conversely, let be given the functions $u_0\in C(S)$ and $u_1\in C(S)$ which satisfy the universal boundary value equations (2), (3). Then there exists the unique function  $u\in C^2(G )\bigcap C^{1,{\rm norm}}(\bar G )$ such that} $$\Delta u=0,\, x\in G$$ $$u|_S=u_0,\quad {{\d u}\over {\d \nu }}|_S=u_1.$$

Here
$$
{{\cos \phi _{x,y}}\over {|x-y|^2}}={{\d }\over {\d \nu
_y}}{1\over {|x-y|}}.
$$

The additional condition of strong convexity on the  surface
$S \in C^3$:

There exists a constant $c_0>0$ such that
$$
|(\nu (u,v),d^2{\bf r}(u,v))|\geq c_0((du)^2+(dv)^2)
$$
for any region
$$
\{ {\bf r}={\bf r}(u,v),\ (u,v)\in D\} =S_1\subset S.
$$
The equations (2), (3) are {\it the universal boundary value equations for Laplace's equation.}

The second part of the theorem is the statement of correct solvability of the universal boundary value problem to the Laplace equation in the space  $C^2(G )\bigcap C^{1,{\rm norm}}(\bar G )$

{\bf Sketch of the proof. Derivation of universal boundary equations.}

Assume that the function $u$ satisfies the conditions:
$$u\in C^2(G)\bigcap C^{1}(\bar G),\,\,\,\Delta u=0,\, x\in G$$.
Then according to the Green formula
$$
u(x)+ {1\over {4\pi }}\int\limits_S{{cos \phi _{xy}}\over {|x-y|^2}}u_0(y)dS_y-{1\over {4\pi }}\int\limits_S{1\over {|x-y|}}u_1(y)dS_y =0,\, x\in G.
$$
By passage to the limit as $x\to S,\, x\in G$ in the Green formula and using the jump formula for the double layer potential
$$
{1\over {4\pi }}\int\limits_S{{cos \phi _{xy}}\over {|x-y|^2}}u_0(y)dS_y \, \to \, {1\over {4\pi }}\int\limits_S{{cos \phi _{xy}}\over {|x-y|^2}}u_0(y)dS_y -{1\over 2}u_0(x)
$$
we obtain the equality (2):
$$
u(x)+ {1\over {2\pi }}\int\limits_S{{cos \phi _{xy}}\over {|x-y|^2}}u_0(y)dS_y-{1\over {2\pi }}\int\limits_S{1\over {|x-y|}}u_1(y)dS_y =0,\, x\in S.
$$

The equality (3) is the consequence of the Ostrogradskii-Gauss theorem.

\bigskip

Now let the function $u$ satisfies the conditions of the theorem:\\
$u\in C^2(G)\bigcap C^{1,norm}(\bar G) $ and $\Delta u=0,\, x\in G$.

Let $S_h$ is a surface which is parallel to the surface $S$, i.e. $S_h$ is the result of displacement of any point of the surface $S$ in the direction of interior normal on the distance $h>0$. Then $S_h$ is $C^2$ surface for any sufficiently small $h$ according to the assumption of the theorem. Let $G_h$ is the subdomain of domain $G$ with the boundary $S_h$. Hence the function $u|_{G_h}\in  C^2(G_h)\bigcap C^{1}(\bar G_h) $ and therefore the functions $u|_{S_h}$ and ${{\d u}\over {\d \nu }}|_{S_h}$ satisfy the equations (2), (3) on the surface $S_h$.

The assumptions of the theorem is sufficient to prove the equality (2), (3) for the functions $u|_{S}$ and ${{\d u}\over {\d \nu }}|_{S}$ by passage to the limit as $h\to +0$.

{\bf Solution of the universal boundary value problem.}

Let the functions $u_0,u_1\in C(S)$ satisfy the universal boundary value equations (2) and (3).

Let us consider the following Neumann problem:
$$
\Delta v=  0\\,\,\, {{\d v}\over {\d \nu }}|_S=u_1
$$
which has the unique solution $v$ according to the condition (3) and  the Lyapunov-Steklov theorem.
Let us set $\hat u_0=v|_S$. Then according to the first statement of the theorem the functions $\hat u_0$ and $u_1$ satisfy the equation (2). The functions $u_0$ and $u_1$ also satisfy this equation according to the assumption of the theorem.

Hence the following equality holds for the function $w=\hat u_0-u_0$:
$$
2\pi w (x)-\int\limits_S w (y){{cos \phi _{xy}}\over {|x-y|^2}}dS_y=0,\, x\in S;$$
then according to the potential theory (\cite{Vladimirov})
$$  w=C_0=const.
$$
Therefore  the function  $u=v-C_0$ is the solution of universal boundary value problem: $$\Delta u=0,\ u|_S=u_0,\ {{\d u}\over {\d \nu }}|_S=u_1.$$

{\bf Remark. } The universal boundary value equations for   Poisson's equation $$\Delta u=f,\ x\in G$$ have the form
\begin{equation}\label {Poisson}
2\pi u_0(x)+\int\limits_S[{{cos\phi _{xy}}
\over {|x-y|^2}}u_0(y)-{{u_1(y)}
\over {|x-y|}}]ds_y-\int\limits_G{{f(y)}
\over {|x-y|}}dy=0,\  x\in S,
\end{equation}
$$
\int\limits_Su_1ds=\int\limits_G{f(y)}dy.
$$

\section { Universal boundary value problem for  the heat equation}

In this section we obtain the universal boundary value equations for the heat equation.

{\bf Theorem 2}.
Let $G= R_+\times  R_+$ be the quarter-plane $\{ (t,x): t,x >0 \}$;
let the function $u=u(t,x)\in C(\bar R_+,C^1(\bar R_+ )
\bigcap C^2(R_+))\bigcap C^1(R_+,C(R_+))$
is bounded on the domain $G$ and satisfies the heat equation:
\begin{equation}\label {heat}
u'_t-u''_{xx}=0,\ (t,x)\in G.x
\end{equation}
Then the function $u$ has the boundary values

$$ v(x)=u(+0,x),\, x>0, $$
$$ \phi (t)=u(t,+0),\, t>0 ,$$
$$\psi (t)=u'_x(t,+0),\, t >0,$$

the boundary values are continuous functions $v,\phi,\psi \in C(R_+)$ and they satisfy the relation
\begin{equation}\label {heatbound}
\phi (t)={1\over {\sqrt {\pi t}}}\int\limits_0^{\infty }e^{-{{\xi ^2}\over {4t}}}v(\xi )d\xi -{1\over {\sqrt {\pi }}}\int\limits_0^{t}{1\over {\sqrt {\tau }}}\psi (t-\tau )d\tau
\end{equation}

The relation (\ref{heatbound}) is the {\it universal boundary value equation for the heat equation on the  quarter-plane $G$}.

\section  {On general boundary value problems}

Here we discuss an interpretation of the universal boundary value equations by using a functional analysis framework.

A rather general boundary value problem for elliptic equations is discussed by Vishik \cite{VMO}  (see also \cite{Herm}) by investigating of extensions of the minimal differential operator or restriction of the maximal differential operator in a Hilbert space.
Different boundary conditions were related in \cite{VMO} with different extensions of the minimal differential operator (self-adjoint extensions, solvable extensions or extensions with another properties), but the universal boundary value equations were not mentioned.

Let $\bf L$ be a linear differential operator of the second order with the domain $D({\bf L})=C^{2}(G )\bigcap C^{1}({\bar G })$. Let the operators ${\bf \gamma }_0,\, {\bf \gamma }_1$ are the trace maps defined on the domain $D({\bf \gamma }_0)=D({\bf \gamma }_1)=D({\bf L})$ by  ${\bf \gamma }_0u=u|_S$, ${\bf \gamma }_1u={{\d u}\over {\d \nu }}|_S$.

Let us introduce an operator ${\cal T}:\ D({\bf L})\, \to \, C(S )
\oplus C(S)\oplus C(G )\equiv Y$
(where $Y$ is the direct sum of Banach spaces) on the domain
$D({\cal T})=D({\bf L})$ by the formula ${\cal T}u=({\bf \gamma }_{0}u,{\bf \gamma }_{1}u,{\bf L }u)$. Instead of the spaces of continuous functions one can use the Sobolev spaces.

The aim of the investigation in \cite{VMO} was to obtain  a description of linear subspaces
$
{\cal P}
$
in the space of boundary values such that the restriction $\bf {L_{{\cal P}}}$ of maximal operator ${\bf L}_{max}$ (without any conditions to boundary values) on the subspace
$\{ u\in D({\bf L}_{max}):\ (\gamma _0u,\gamma _1u)\in {\cal P}\}$ has the bounded inverse operator or compact inverse operator. The general boundary value problem according to \cite{VMO} is the elliptic partial differential equation together with the equation on boundary values which specify the subspace $\cal P$.

The problem of finding the solvable restrictions of the maximal elliptic operator ${\bf L}_{max}$  is reduced to the problem of describing  the set of operators ${\bf A}, \, {\bf B}$ in the space $C(S)$ such that the intersection of three subspaces in the space $Y$: ${\rm Im}({\cal T})$,
$\{ (v_0,v_1,f)\in Y:\ {\bf A}v_{0}+{\bf B}v_{1}=\theta \}$ and $\{ (v_0,v_1,f)\in Y:\ f=\hat f\}$ consists of the unique point for arbitrary $\hat f \in C(G)$. The third component of the space $Y$ has the special role in the approach of \cite{VMO}. The aim of our investigation is to present the image ${\rm Im}({\cal T})$ by the linear equation on the space $Y$  (the equation (7) below) such that all of three components of the image of operator $\cal T$ are having the same rights in this equation.

Our investigation of the universal boundary value problem  for Laplace's equation (1) and  for the heat equation (4) in principle can be extended to the universal boundary value problem for more general linear differential equations. One can proceed as follows.
The universal boundary equation for the differential equation $${\bf L}u=f$$
can be defined as an equation  in the space $Y$:
\begin{equation}\label{Y}
\{ (v_0,v_1,f)\in Y:\ {\bf A}v_{0}+{\bf B}v_{1}+{\bf C}f=\theta \},
\end{equation}
where linear operators ${\bf A},{\bf B},{\bf C }$ are defined on the spaces $C(S )$, $ C(S )$, $ C(G )$ respectively and take values in some Banach space $Z$. Of course, the main task is to find an explicite form of the
operators ${\bf A},{\bf B},{\bf C }$.

For example the equation (7) for Poisson equation $\Delta u=f,\ x\in G$ is the generalization of equations (2), (3) in the form of equation (7):
$$
2\pi u_0(x)+\int\limits_S[{{cos\phi _{xy}}\over {|x-y|^2}}u_0(y)-{{u_1(y)}
\over {|x-y|}}]ds_y-\int\limits_G{{f(y)}
\over {|x-y|}}dy=0,\  x\in S.
$$
$$
\int\limits_Su_1ds=\int\limits_G{f(y)}dy.
$$

\bigskip

The universal boundary value problem for the homogeneous equation  is the problem of description of the kernel of the maximal operator ${\bf L}$ in terms of boundary values of the elements of this kernel,
\begin{equation}\label{hom}
u\in {\rm Ker}({\bf L}) \ \Leftrightarrow \ \{ (v_0,v_1)\in C(S)\times C(S):\ {\bf A}v_0+{\bf B}v_1=0\}.
\end{equation}
Eq (\ref{hom}) is an abstract form of the universal boundary value equations (2), (3) for Laplace equation (1).

In \cite{Burs} the problem of describing relations between boundary values is discussed and by using the Fourier transform
an infinite number of equations for some functionals from  the boundary functions is obtained.

\section  {Discussions and  Conclusions}

In this note Laplace's and the heat equations are considered and the  restrictions on  the boundary values  of the solutions  
and its normal derivatives are studied.  The linear integral equations for the boundary values of the solution
and its normal derivatives are obtained   which are called the universal boundary value equations.

The universal boundary value problem appeared in the cosmological considerations \cite{ArVol}. One had to define an operator
$e^{\tau \frac{d^2}{dt^2}}v(t)=u(\tau,t)$ without specifying boundary conditions on the quarter-plane $\tau, t >0$. To this end one of the authors (I.V.)  obtained the universal boundary value equations for the heat and Laplace's equations. Concerning applications of the  operator $e^{\tau \frac{d^2}{dt^2}}v(t)$ on the half-plane to nonlinear equations see \cite{Vladimirov-2}.

Restrictions on the Cauchy data are mentioned in the studying the wave equation on non-globally-hyperbolic manifold \cite{GGKV}. In this work the description of the kernel (8) for operator of Cauchy problem for wave equation on  non-globally-hyperbolic manifold is obtained.

There could be  important applications of the universal boundary value problem in quantum mechanics when the Schrodinger equation with the degenerated Hamiltonian is considered \cite{Sakbaev-1, Sakbaev-2}. In this paper it was obtained a description of the kernel (\ref{hom}) for operator of the Cauchy problem for the Schrodinger equation with degenerated Hamiltonian (see theorem 2.2 in \cite{Sakbaev-2}) and solution of an optimal problem (see theorem 11.1 in \cite{Sakbaev-2}) on the set (\ref{Y}) for this equation.

It would be interesting to obtain generalizations of
the universal boundary value equations considered in this paper to other partial differential equations.

\section{Acknowledgements}

We are grateful to I.Ya. Aref'eva, Yu.N. Drozzinov, L.S. Efremova, A.K. Gushchin, V.P. Mikhailov, N.N. Shamarov, V.S. Vladimirov and V.V. Zharinov for useful discussions. The paper was partially supported by the grant NSh-2928.2012.1.

\end{document}